\documentclass[a4paper, 12pt, twoside]{article}
\usepackage{amsmath}
\usepackage{amsthm}
\usepackage{amsfonts}
\usepackage{pifont}
\usepackage{color}
\usepackage{bbm}
\usepackage{mathrsfs}
\usepackage{fancyhdr}
\usepackage{graphicx}

\usepackage{psfrag}
\usepackage{time}
\usepackage{txfonts}

\usepackage{stmaryrd}
\usepackage{mathrsfs}
\usepackage{txfonts}
\usepackage{amsfonts}

\usepackage{indentfirst}

\newcommand{\xiaosan}{\fontsize{15pt}{22pt}\selectfont}
\newcommand{\sihao}{\fontsize{14pt}{21pt}\selectfont}

\newcommand{\xiaosi}{\fontsize{12pt}{18pt}\selectfont}

\usepackage{curves}

\numberwithin{equation}{section}

\newtheorem{theorem}{ {Theorem}}[section]

\newtheorem{remark} {   {Remark}}[section]
\newtheorem{corollary} {  {Corollary}}[section]
\newtheorem{lemma} {  {Lemma}}[section]
\newtheorem{definition}{  {Definition}}[section]

\begin{document}
\setlength{\parindent}{2em}
\newpage
\fontsize{12}{22}\selectfont\thispagestyle{empty}
\renewcommand{\headrulewidth}{0pt}
 \lhead{}\chead{}\rhead{} \lfoot{}\cfoot{}\rfoot{}
\noindent

\title{\bf Anderson  Localization   for  the Almost Mathieu Operator  in Exponential Regime}
\author{{Wencai  Liu and Xiaoping Yuan*}\\
{\em\small School of Mathematical Sciences}\\
{\em\small Fudan University}\\
{\em\small  Shanghai 200433, People's Republic of China}\\
{\small 12110180063@fudan.edu.cn}\\
{\small *Corresponding author: xpyuan@fudan.edu.cn}}
\date{}
\maketitle

\renewcommand{\baselinestretch}{1.2}
\large\normalsize

\begin{abstract}
 For the
almost Mathieu operator $(H_{\lambda,\alpha,\theta}u)_n=u_{n+1}+u_{n-1}+2\lambda \cos2\pi(\theta+n\alpha)u_n$,
 Avila and  Jitomirskaya guess that for a.e.\;$\theta$, $H_{\lambda,\alpha,\theta}$ satisfies Anderson localization  if $ |\lambda| >  e^{  \beta} $,
  and they   establish  this for $ |\lambda| >  e^{\frac{16}{9} \beta}$.
  In the present paper, we extend    their result
  to    regime $ |\lambda| >  e^{\frac{3}{2} \beta}$.
\end{abstract}

\setcounter{page}{1} \pagenumbering{arabic}\topskip -0.82in
\fancyhead[LE]{\footnotesize  Introduction}
\section{\xiaosan \textbf{Introduction}}
    The almost Mathieu operator (AMO) is the (discrete) quasi-periodic   Schr\"{o}dinger operator  on  $   \ell^2(\mathbb{Z})$:
 \begin{equation}\label{G11}
 (H_{\lambda,\alpha,\theta}u)_n=u_{n+1}+u_{n-1}+ \lambda v(\theta+n\alpha)u_n,   \text{ with }  v(\theta)=2\cos2\pi \theta,
 \end{equation}
where $\lambda$ is the coupling, $\alpha $ is the frequency, and $\theta $ is the phase.
\par
  $ H_{\lambda,\alpha,\theta}$ is a tight binding model for the Hamiltonian of an electron in a one-dimensional
lattice or
in a two-dimensional lattice, subjecting to a perpendicular (uniform) magnetic field (through a Landau gauge)$ {\cite{Har},\cite{Rau}}$.
This model   also describes a square lattice with anisotropic nearest neighbor coupling and
isotropic next nearest neighbor coupling, or anisotropic coupling to the nearest neighbors and next nearest neighbors
on a triangular lattice $ {\cite{Bel},\cite{Tho}}$.  For more applications in physics,  we refer the reader to
 \cite{Las}  and the references therein.
 \par
Besides  its relations to some fundamental problems in physics, the AMO itself is also
fascinating  because of its remarkable richness of the related spectral theory. In Barry Simon's list of Schr\"{o}dinger operator problems
for the twenty-first century $   \cite{Sim} $, there are three problems   about the AMO. The spectral theory of AMO
has attracted  many   authors, for example,  Avila-Jitomirskaya\cite{AJ1}, \cite{AJ2},  Avila-Krikorian\cite{AK},  Bourgain\cite{B1},\cite{B2},  Jitomirskaya-Simon \cite{JS} and so on.
\par
Anderson localization (i.e., only pure point spectrum with exponentially decaying eigenfunctions) is not only
meaningful in physics, but also relates to some problems of   the quasi-periodic   Schr\"{o}dinger operator, such as the
reducibility of cocycles   via   Aubry duality    \cite{GJLS}  and
the  Ten Martini Problem (Cantor spectrum conjecture)   \cite{AJ1}.
\par
For $\alpha\in \mathbb{Q}$,  it is easy to verify that $ H_{\lambda,\alpha,\theta}$ has no  eigenvalues,
let alone Anderson localization. Thus, in the present paper,  we always assume  $\alpha\in \mathbb{R}\backslash \mathbb{Q}$.
\par
For simplicity, we say $H_{\lambda,\alpha,\theta}$ satisfies AL if  for  a.e.\;phase $\theta$, $H_{\lambda,\alpha,\theta}$  satisfies Anderson localization.
\par
Avila and Jitomirskaya guess that    $H_{\lambda,\alpha,\theta}$ satisfies AL for $ |\lambda| >  e^{  \beta}$ (Remark 9.2,  $\cite{AJ1}$),
where
 \begin{equation}\label{G12}
  \beta= \beta(\alpha)=\limsup_{n\rightarrow\infty}\frac{\ln q_{n+1}}{q_n},
 \end{equation}
 and  $ \frac{p_n}{q_n} $  is  the continued fraction approximants   to $\alpha$.
  One usually calls set $\{\alpha \in \mathbb{R}\backslash \mathbb{Q}|\;\beta(\alpha)>0\}$    exponential regime and set $\{\alpha\in \mathbb{R}\backslash \mathbb{Q}|\;\beta(\alpha)=0\}$   sub-exponential regime.
\par
This guess is optimal in some way. On the one hand, for every $\alpha$ there is a generic set of $\theta$ for which there
is  no  eigenvalues \cite{JS}. On the other hand,  if $|\lambda|\leq e^{\beta}$, for every $ \theta$, $ H_{\lambda,\alpha,\theta}$ has no localized eigenfunctions (i.e., exponentially decaying eigenfunctions) \cite{G}.
\par
  In $ \cite{BJ2}$,  Bourgain and Jitomirskaya  prove that   $H_{\lambda,\alpha,\theta}$
  satisfies AL if
$ \alpha\in DC$\footnote{We say $\alpha \in \mathbb{R}\backslash \mathbb{Q} $ satisfies a Diophantine condition $\text{DC}(\kappa,\tau)$ with $\kappa>0$ and $\tau>0$,
if
$$ |q\alpha-p|>\kappa |q|^{-\tau}  \text{ for   any } (p,q)\in \mathbb{Z}^2, q\neq 0.$$
Let $\text{DC}=\cup_{\kappa>0,\tau>0}\text{DC}( \kappa,\tau)$. We say $\alpha$  satisfies   Diophantine condition, if $\alpha\in \text{DC}$.
 Notice that  $\beta(\alpha)=0$  for
$\alpha\in DC$.} and $|\lambda|>1$. Avila and Jitomirskaya     obtain that
    $H_{\lambda,\alpha,\theta}$
  satisfies AL if  $\beta(\alpha)=0$ and $|\lambda|>1$ \cite{AJ2}. In fact,  Avila and Jitomirskaya's  analysis
 also suggests that    $H_{\lambda,\alpha,\theta}$  satisfies  AL if $ |\lambda| >  e^{C \beta }$,
 where $C$ is a large  absolute constant
(after carefully checking their proof ).
  In  $ \cite{AJ1}$, Avila and Jitomirskaya give a definite quantitative description of the constant $C$  and get  $C=\frac{16}{9}$.
 In the present paper, we extend  to
   regime $ |\lambda| >  e^{   \frac{3\beta}{2}}$, i.e., the following theorem.

 \begin{theorem} (\textbf{Main Theorem})\label{Main Theorem}
 Let $ \alpha\in \mathbb{R}\backslash\mathbb{Q}$ be such that $ \beta=\beta(\alpha)<\infty$,
   then for
  almost every phase $\theta $, $ H_{\lambda,\alpha,\theta}$ satisfies  Anderson localization   if $| \lambda |>  e^{\frac{3}{2} \beta}$.
\end{theorem}
 Here we would like to talk about some  histories of the  investigation  to  Anderson localization in more details.
 To state the problem more simply,  we sometimes  drop the parameters dependence, such as $\lambda,\alpha,\theta$ and so on.
 \par
 Let $H=H_{\lambda,\alpha,\theta}$.
  Define $H_I=R_I{H}R_I$, where $R_I=$ coordinate restriction to $I=[x_1,x_2]\subset\mathbb{Z}$, and denote by
$ {G}_I =(  {H}_I-E)^{-1}$  the associated Green  function, if   $ {H}_I-E$  is invertible. Denote by $ {G}_{I}(x,y)$    the matrix elements of  Green  function
${G}_I$. Note that ${G}_I$ depends on $\lambda,\alpha,\theta,E$.
\par
It is  easy to check  if the Green function $G_I(\theta)$ satisfies
\begin{equation}\label{G13}
  |{G}_{I}(\theta)(m,n)|<e^{-c|m-n|} \text { for } |m-n|> |I|/ 5,
\end{equation}
where $c>0$ and $|I|=b-a+1$  for $I=[a,b]$,
then Anderson localization  holds.
Unfortunately, (\ref{G13}) does not hold in general.

 Nevertheless,  Bourgain  proves that (\ref{G13}) holds for $I=[0,N]$      except for $\theta$ in a small exceptional set. A typical
  statement would be the following
  \begin{equation}\label{G14}
    ||G_{[0,N]}(\theta)||<N^{1-\delta}
  \end{equation}
  and
  \begin{equation}\label{G15}
     |{G}_{[0,N]}(\theta)(m,n)|<e^{-c|m-n|} \text { if } |m-n|> N/ 5
  \end{equation}
  for all $\theta$ outside a set of measure $<e^{-N^\sigma}$ if $|\lambda|>1$.
  Here $ \delta, \sigma$ are some positive constants.
 Via Bourgain's careful  arguments,  he proves that  for  a full
Lebesgue measure subset of Diophantine frequencies, $H_{\lambda,\alpha,\theta}$
satisfies AL if $|\lambda|>1$.  See Bourgain's book  \cite{B2}  for details.
 \par
 In  \cite{BJ2}, Bourgain and Jitomirskaya develop another  subtle  way to  set up sharp estimate of Green function. We recall the main idea.
 For any $k>0$, they  success  to look for a interval $I=[x_1,x_2]\subset \mathbb{Z}$
  with $k\in I$ and $\text{dist}(k,x_i)>|I|/ 5$,  such that
 \begin{equation}\label{G16}
    |G_I(x_i,k)|<e^{-c|k-x_i|} \text{ for some } c>0.
 \end{equation}
 Then Anderson localization  follows from (\ref{G16}) in a well known manner--block resolvent expansion (see \cite{B1} for example).
 As a result, they display  AL for $H_{\lambda,\alpha,\theta}$
   if
$ \alpha\in DC$ and $|\lambda|>1$.
 Their  discussion  strongly
relies   on the cosine potential. Concretely,  their methods can only apply to
quasi-periodic   Schr\"{o}dinger operator  (\ref{G11}) with $v=2\cos2\pi \theta$.
 How to apply to general potential $v$ is still open.
\par
Following the program of  Bourgain-Jitomirskaya in \cite{BJ2},  Avila and Jitomirskaya estimate the Green function more    finely
 \cite{AJ2}.
In addition using  Lemma \ref{Le25} below technically, Avila and Jitomirskaya
   obtain that
    $H_{\lambda,\alpha,\theta}$
  satisfies AL for   $\beta(\alpha)=0$ and $|\lambda|>1$.
  Furthermore,
in another paper\cite{AJ1}, they distinguish $k$  resonance and   non-resonance respectively to  look for interval $I$ such that (\ref{G16}) holds. Together with  some results  in  \cite{AJ2},\cite{BJ2},
they prove that AL holds if $ |\lambda| >  e^{   \frac{16\beta}{9}}$.
\par
We investigate  the Anderson localization as the program of Avila and Jitomirskaya in \cite{AJ1}. If  $k$ is non-resonant,  Avila and Jitomirskaya's  analysis is optimal, thus
we use directly (Theorem \ref{Th24}).  In the present paper,  we  focus   our attention   on  the     resonant $k$,  and  carry on   more
 subtle computation in estimating   Green function.
\par
 The present paper is organized as follows:
 \par
In \S 2,  we give some preliminary notions and facts which are taken from other authors, such as  Avila-Jitomirskaya  \cite{AJ1},  Bourgain\cite{B2} and so on. In \S3, we set up the regularity of  resonant $y$ if  $ |\lambda| >  e^{   \frac{3\beta}{2}}$.
 In \S4, we give the proof of Main theorem by block resolvent expansion.
\section{ Preliminaries and some known results  }
 \par
  It is well known that  Anderson localization for a self-adjoint operator  $H$ on $\ell^2(\mathbb{Z})$ is equivalent to  the following statements.
  \par
  Assume $\phi$ is an extended state, i.e.,
  \begin{equation}\label{G21}
    H\phi=E\phi  \text{ with  } E\in \Sigma(H) \text{ and }  |\phi(k)|\leq (1+|k|)^C,
  \end{equation}
  where $\Sigma(H)$ is the spectrum of   self-adjoint   operator $H$.
  Then there exists some constant $c>0$ such that
  \begin{equation}\label{G22}
    | \phi(k)|< e^{-c|k|} \text{ for } k\rightarrow \infty.
  \end{equation}
 \par
 The above statements can be proved by Gelfand-Maurin Theorem.
 See \cite{BER} for the proof of    continuous-time Schr\"{o}dinger operator.
The proof of  discrete     Schr\"{o}dinger operator is similar, see   \cite{LIU2} for example.

We will actually prove a slightly more precise version of Theorem $ \ref{Main Theorem}$. Let
\begin{equation}\label{G23}
  \mathscr{R}_1=\{ \theta:|\sin \pi(2\theta+k \alpha )|\leq k^{-2} \text{ holds  for  infinitely  many } k, k\in \mathbb{Z}\},
\end{equation}
and $\mathscr{R}_2= \{\theta: \exists s \in \mathbb{Z} \text{ such  that } 2\theta +s\alpha \in \mathbb{Z}\}$.
Clearly,  $  \mathscr{R}= \mathscr{R}_1\cup \mathscr{R}_2$  has zero Lebesgue measure.
\begin{theorem}\label{Th21}
Let $ \alpha\in \mathbb{R}\backslash\mathbb{Q}$ be such that $ \beta=\beta(\alpha)<\infty$,
then   $ H_{\lambda,\alpha,\theta}$  satisfies   Anderson localization if  $\theta \notin \mathscr{R}$ and  $|\lambda|>e^{\frac{3\beta}{2}}$.
\end{theorem}

 \par
If   $ \alpha $ satisfies  $ \beta(\alpha) =0$,   Theorem   \ref{Th21}  has been proved by  Avila-Jitomirskaya  in $\cite{AJ1}$ and $\cite{AJ2} $. Thus in the present  paper, we fix $ \alpha\in \mathbb{R}\backslash\mathbb{Q}$
  such that $ 0<\beta(\alpha) <\infty$.  Unless  stated otherwise,   we always  assume $\lambda>e^{\frac{3}{2}\beta}$ ( for $\lambda<-e^{\frac{3}{2}\beta}$,  notice that   $H_{\lambda,\alpha,\theta}=H_{-\lambda,\alpha,\theta+\frac{1}{2}}$), $\theta\notin  \mathscr{R}$
and $E\in \Sigma_{\lambda,\alpha}$\footnote{ The spectrum  of operator
 $ H_{\lambda,\alpha,\theta}  $
 does not depend on $ \theta$,  denoted by  $  \Sigma_{\lambda,\alpha}$. Indeed, shift is an unitary
 operator on $\ell^2(\mathbb{Z})$, thus $  \Sigma  _{\lambda,\alpha,\theta}= \Sigma  _{\lambda,\alpha,\theta+\alpha}$, where $  \Sigma _{\lambda,\alpha,\theta}$
 is the spectrum of  $ H_{\lambda,\alpha,\theta}  $. By  the minimality of $ \theta\mapsto \theta+\alpha$ and  continuity of  spectrum
  $  \Sigma  _{\lambda,\alpha,\theta}$ with respect  to  $\theta$, the statement follows.}. Since this does not   change  any of the statements,  sometimes the dependence of parameters $E,\lambda,\alpha, \theta$ will be ignored in the following.
  \par
  Given      an extended state $\phi$ of $H_{\lambda,\alpha,\theta}$, without loss of generality assume $\phi(0)=1$.
  Our objective is to prove that there exists some $c>0$ such that
  $$| \phi(k)|< e^{-c|k|} \text{ for } k\rightarrow \infty.$$

\par
Let us denote
$$ P_k(\theta)=\det(R_{[0,k-1]}(H_{\lambda,\alpha,\theta}-E) R_{[0,k-1]}).$$
Following $ {\cite{JKS}}$, $P_k(\theta)$ is an even function of $ \theta+\frac{1}{2}(k-1)\alpha$  and can be written as a polynomial
      of degree $k$ in $\cos2\pi (\theta+\frac{1}{2}(k-1)\alpha )$ :

     \begin{equation}\label{G24}
      P_k(\theta)=\sum _{j=0}^{k}c_j\cos^j2\pi (\theta+\frac{1}{2}(k-1)\alpha)    \triangleq  Q_k(\cos2\pi  (\theta+\frac{1}{2}(k-1)\alpha)).
     \end{equation}

     Let $A_{k,r}=\{\theta\in\mathbb{R} \;|\;Q_k(\cos2\pi   \theta   )|\leq e^{(k+1)r}\} $ with $k\in \mathbb{N}$ and $r>0$.
\begin{lemma}$ (\text{p}. 16, \cite{AJ1})$ \label{Le22}
The following inequality  holds
    \begin{equation}\label{G25}
     \lim_{k\rightarrow\infty}\sup_{\theta\in \mathbb{R}} \frac{1}{k} \ln | P_k(\theta)|\leq \ln \lambda.
    \end{equation}
\end{lemma}

    \par
    By Cramer's rule (p. 15, $\cite{B2}$)  for given  $x_1$ and $x_2=x_1+k-1$, with
     $ y\in I=[x_1,x_2] \subset \mathbb{Z}$,  one has
     \begin{eqnarray}
       |G_I(x_1,y)| &=&  \left| \frac{P_{x_2-y}(\theta+(y+1)\alpha)}{P_{k}(\theta+x_1\alpha)}\right|,\label{G26}\\
       |G_I(y,x_2)| &=&\left|\frac{P_{y-x_1}(\theta+x_1\alpha)}{P_{k}(\theta+x_1\alpha)} \right|.\label{G27}
     \end{eqnarray}
By Lemma \ref{Le22}, the numerators in  (\ref{G26}) and (\ref{G27}) can be bounded uniformly with respect to $\theta$. Namely,
for any $\varepsilon>0$,
\begin{equation}\label{G28}
    | P_n(\theta)|\leq e^{(\ln \lambda+\varepsilon)n}
\end{equation}
for $n$ large enough.
\begin{definition}\label{Def21}
Fix $t > 0$. A point $y\in\mathbb{Z}$ will be called $(t,k)$-regular if there exists an
interval $[x_1,x_2]$  containing $y$, where $x_2=x_1+k-1$, such that
\begin{equation}\label{G29}
  | G_{[x_1,x_2]}(y,x_i)|<e^{-t|y-x_i|} \text{ and } |y-x_i|\geq \frac{1}{5} k \text{ for }i=1,2;
\end{equation}
otherwise, $y$ will be called $(t,k)$-singular.
\end{definition}
It is  easy to check that (p. 61, $\cite{B2}$)
 \begin{equation}\label{G210}
   \phi(x)= -G_{[x_1 ,x_2]}(x_1,x ) \phi(x_1-1)-G_{[x_1 ,x_2]}(x,x_2) \phi(x_2+1),
 \end{equation}
 where  $ x\in I=[x_1,x_2] \subset \mathbb{Z}$.
Our strategy is to establish the $(t,  k(y) )$-regular for every large $y$,  then localized property is easy to
obtain by   $(\ref{G210})$ and the block resolvent expansion.
       \begin{definition}
     We  say that the set $\{\theta_1, \cdots ,\theta_{k+1}\}$ is $ \epsilon$-uniform if
      \begin{equation}\label{G211}
        \max_{ x\in[-1,1]}\max_{i=1,\cdots,k+1}\prod_{ j=1 , j\neq i }^{k+1}\frac{|x-\cos2\pi\theta_j|}
        {|\cos2\pi\theta_i-\cos2\pi\theta_j|}<e^{k\epsilon}.
      \end{equation}
     \end{definition}
      \begin{lemma}\label{Le23}$(\text{Lemma 9.3 },\cite{AJ1})$
      Suppose  $\{\theta_1, \cdots ,\theta_{k+1}\}$ is  $ \epsilon_1$-uniform. Then there exists some $\theta_i$ in set  $\{\theta_1, \cdots ,\theta_{k+1}\}$ such that
     $\theta_i\notin A_{k,\ln\lambda-\epsilon}$ if    $ \epsilon>\epsilon_1$ and $ k$
      is sufficiently large.
      \end{lemma}
   Assume without loss of generality that $y>0$. Define $b_n=q_n^{8/9}$, where $q_n$ is given by (\ref{G12}), and  find $n$ such that $b_n\leq y<b_{n+1}$.
        We will distinguish two cases:
         \par
        (i)   $|y-\ell q_n|\leq b_n$ for some $\ell\geq1$,  called resonance.
          \par
        (ii)     $|y-\ell q_n|> b_n$ for all $\ell\geq0$, called  non-resonance.
\par
For the  non-resonant $y$,  Avila and Jitomirskaya have established  the regularity for $y$, which is optimal.  We  give the theorem  directly.

\begin{theorem}$(\text{Lemma 9.4}, \cite{AJ1})$\label{Th24}
Assume $\theta \notin \mathscr{R}$,  $  \lambda>e^{\beta }$ and   $y$ is non-resonant. Let $s\in\mathbb{N}$   be
the largest number such that $ sq_{n-1}\leq dist(y,\{ \ell q_n\}_{\ell\geq0}) $,
then $\forall \varepsilon>0$,
$ y$ is $  (\ln\lambda+9\ln (s q_{n-1}/q_n)/q_{n-1}-\varepsilon,2sq_{n-1}-1)$-regular if  $ n$ is  large enough (or equivalently $y$    is  large enough). In particular,
$ y$ is $  (\ln\lambda-\beta-\varepsilon,2sq_{n-1}-1)$-regular.
\end{theorem}

    \begin{lemma}$(\text{Lemma }9.8, \cite{AJ1})$\label{Le25}
Let $m \in  \mathbb{N}$ be such that
$m<\frac {q_{r+1}} {10 q_n}$, where $r\geq n$.
Given a integer  sequence $|m_k|\le m-1,$ $\,k=1,\cdots,q_n,$ let $1\le k_0\le q_n$  be such that

\begin{equation}\label{G212}
|\sin  \pi(x+ (k_0+ m_{k_0} q_r)\alpha   )|=
\min_{1\le k\le q_n}|\sin  \pi(x+ (k+m_k q_r)\alpha )|,
\end{equation}
then
\begin{equation}\label{G213}
\left |\sum_{^{k=1}_{k\not=k_0}}^{q_n}\ln|\sin
   \pi(x+ (k+m_k q_r)\alpha  )|+(q_n-1)\ln 2 \right | <
C\ln q_n+C(\Delta_n+(m-1)\Delta_r) q_n\ln q_n,
\end{equation}
where  $\Delta_n= |q_n\alpha-p_n|$.
\end{lemma}

\section{Regularity for   resonant $y$  }
In this section, we mainly concern the regularity  for   resonant $y$.
If $b_n\leq y<b_{n+1}$ is  resonant, by the definition of resonance, there exists some  positive integer $ \ell$ with   $1 \leq \ell\leq q_{n+1}^{8/9}/q_n$
such that   $|y-\ell q_n|\leq b_n$.  Fix the positive integer $\ell$   and set $I_1, I_2\subset \mathbb{Z}$ as follows
\begin{eqnarray*}
  I_1 &=& [-[\frac{2}{3}q_n],[\frac{2}{3}q_n]-2], \\
   I_2 &=& [(\ell-1) q_n+[\frac{2}{3}q_n]-1,(\ell+1)q_n-[\frac{2}{3}q_n]-1 ],
\end{eqnarray*}
and let $\theta_j=\theta+j\alpha$ for $j\in I_1\cup I_2$. The set $\{\theta_j\}_{j\in I_1\cup I_2}$
consists of $2q_n$ elements.
\par
Note that, below,     we       replace $I=[x_1,x_2]\cap \mathbb{Z}$  with  $I=[x_1,x_2]$ for simplicity, and assume  $\varepsilon>0$  is  sufficiently  small.
     \par
       We will use the following three steps  to establish regularity for $y$.
      \textbf{ Step 1}:  We  set up the $ \frac{\beta}{2}+\varepsilon$-uniformity of $\{\theta_j\}$ where $\theta_j=\theta+j\alpha $
    and $ j$ ranges through $I_1\cup I_2$. By Lemma $ \ref{Le23}$, there exists some $j_0$ with  $j_0\in I_1\cup I_2$
    such that
      $ \theta_{j _0}\notin  A_{2q_n-1,\ln\lambda-\frac{\beta}{2}-2\varepsilon }$.   \textbf{ Step 2}: We   show that
      $\forall j \in I_1, \theta_j \in   A_{2q_n-1, \ln\lambda-\frac{\beta}{2}-2\varepsilon }$ if $\lambda>e^{\frac{3}{2}\beta}$.  Thus there exists
       $ \theta_{j_0} \notin  A_{2q_n-1,\ln\lambda-\frac{\beta}{2}-2\varepsilon}$ for some $j_0 \in I_2$.
        \textbf{ Step 3}:
       We establish the regularity for $y$.
\begin{remark}\label{Remark44}
In \cite{AJ1}, Avila and Jitormirskaya construct $ I_1 =  [-[\frac{5}{8}q_n],[\frac{5}{8}q_n]-1],
   I_2 =  [(\ell-1) q_n+[\frac{5}{8}q_n],(\ell+1)q_n-[\frac{5}{8}q_n]-1 ] $ and set $\theta_j=\theta+j\alpha$ for $j\in I_1\cup I_2$.
They use the above three steps to establish the   regularity of $y$. More precisely, firstly, they
  establish the $ \frac{\beta}{2}+\varepsilon$-uniformity of $\{\theta_j\}$ and there exists
      $ \theta_{j _0}\notin  A_{2q_n-1,\ln\lambda-\frac{\beta}{2}-2\varepsilon }$  for some $j_0\in I_1\cup I_2$.  Secondly,   they    prove that
      $\forall j \in I_1, \theta_j \in   A_{2q_n-1, \ln\lambda-\frac{\beta}{2}-2\varepsilon }$    and  thus there exists
       $ \theta_{j_0} \notin  A_{2q_n-1,\ln\lambda-\frac{\beta}{2}-2\varepsilon}$ for some $j_0 \in I_2$, if $ \lambda >e^{\frac{16}{9}\beta}$.
        Thirdly,
       they set up the regularity of $y$. In the present paper, we reconstruct $I_1$ and $I_2$, and  show that  the three steps
      also hold.
\end{remark}
Recall that
\begin{equation}\label{G31}
\forall 1\leq k <q_{n+1}, \| k\alpha\|_{\mathbb{R}/\mathbb{Z}}\geq  \Delta_n,
\end{equation}

and
\begin{equation}\label{G32}
      \frac{1}{2q_{n+1}}\leq\Delta_n \leq\frac{1}{q_{n+1}},
\end{equation}
where $||x||_{\mathbb{R}/\mathbb{Z}}=\min_{j\in \mathbb{Z}}|x-j|$.
\par
\textbf{Step 1:} We   establish the $(\frac{\beta}{2}+\varepsilon )$-uniformity for $\{\theta_j\}_{j\in I_1 \cup I_2}$.
\par
 In Lemma $\ref{Le25}$, let $r=n$ and  $m=\ell\leq q_{n+1}^{8/9}/q_n$,   one has
 \begin{equation*}
  ( \Delta_n+ (m-1)\Delta_r)q_n=\ell \Delta_n q_n\leq C,
 \end{equation*}
 since $\Delta_n \leq\frac{1}{q_{n+1} }$ by (\ref{G32}).
  Moreover, we obtain  the following
lemma.
\begin{lemma}\label{Le32}
Given a integer sequence $|m_k|\le \ell-1,$ $\,k=1,\cdots,q_n,$ let $1\le k_0\le q_n$  be such that

\begin{equation}\label{G33}
|\sin  \pi(x+ (k_0+ m_{k_0} q_n)\alpha   )|=
\min_{1\le k\le q_n}|\sin  \pi(x+ (k+m_k q_n)\alpha )|,
\end{equation}
then
\begin{equation}\label{G34}
 -(q_n-1)\ln 2 -C\ln q_n\leq \sum_{^{k=1}_{k\not=k_0}}^{q_n}\ln|\sin
   \pi(x+ (k+m_k q_n)\alpha  )| \leq -(q_n-1)\ln 2+
C\ln q_n.
\end{equation}
\end{lemma}
\begin{theorem}   \label{Th33}
 $\forall$ $ \varepsilon >0$,    the set $\{\theta_j\}_{j\in I_1\cup I_2}$
is $(\frac{\beta}{2}+\varepsilon )$-uniform for $\theta \notin  \mathscr{R}  $ and sufficiently large $n$.
\end{theorem}
\textbf{Proof:} Let $$I'_1  =  [-[\frac{2}{3}q_n],-[\frac{2}{3}q_n]+ q_n-1]$$ and $$
   I'_2  =  [-[\frac{2}{3}q_n]+ q_n ,  [\frac{2}{3}q_n]-2 ]\cup [(\ell-1) q_n+[\frac{2}{3}q_n]-1,(\ell+1)q_n-[\frac{2}{3}q_n]-1 ].$$
 Clearly,  both  $\{\theta_j\}_{j\in I'_1}$ and $\{\theta_j\}_{j\in I'_2}$
consist  of $ q_n$ elements, and $ I'_1\cup  I'_2 =I _1\cup  I _2$.
In ($\ref{G211}$), let $x=\cos2\pi a$, $k=2q_n-1$  and take    the logarithm.  Thus in order to prove the theorem,
 it suffices to show that  for any $a\in \mathbb{R}$ and $i\in I'_1\cup  I'_2$,

  $$ \ln \prod_{ j\in I'_1\cup  I'_2  , j\neq i } \frac{|\cos2\pi a-\cos2\pi\theta_j|}
        {|\cos2\pi\theta_i-\cos2\pi\theta_j|}\;\;\;\;\;\;\;\;\;\;\;\;\;\;\;\;\;\;\;\;\;\;\;\;\;\;\;\;\;\;\;\;\;\;\;\;\;\;\;\;\;\;\;\;\;\;\;\;\;\;\;\;\;
        \;\;\;\;\;\;\;\;\;\;\;\;\;\;\;$$
          $$=   \sum _{ j\in I'_1\cup  I'_2  , j\neq i }\ln|\cos2\pi a-\cos2\pi \theta_j|- \sum _{ j\in I'_1\cup  I'_2  , j\neq i }\ln|\cos2\pi\theta_i  -\cos2\pi \theta_j| $$
       \begin{equation}\label{2G35}
       <   (2q_n- 1)(\frac{\beta}{2}+\varepsilon).\;\;\;\;\;\;\;\;\;\;\;\;\;\;\;\;\;\;\;\;\;\;\;\;\;\;\;\;\;\;\;\;\;\;\;\;\;\;\;\;\;\;\;\;\;\;\;\;
     \;\;\;\;\;\;\;\;\;\;\;\;\;\;\;\;\;\;\;\;\;
       \end{equation}

  Without loss of generality assume $i \in I'_1$.
 We  estimate $ \sum _{ j\in I'_1\cup  I'_2  , j\neq i }\ln|\cos2\pi a-\cos2\pi \theta_j| $ first.
\par
Clearly,
   $$ \sum _{ j\in I'_1\cup  I'_2  , j\neq i }\ln|\cos2\pi a-\cos2\pi \theta_j| \;\;\;\;\;\;\;\;\;\;\;\;\;\;\;\;\;\;\;\;\;\;\;\;\;\;\;\;\;\;\;\;\;\;\;\;\;\;\;\;\;\;\;\;\;\;\;\;\;\;\;\;\;\;\;\;$$
$$\;\;\;\;\;\;\;\;\;\;\;\;\;\;\;\;\;\;\;\;\;=\sum_{ j\in I'_1\cup  I'_2  , j\neq i }\ln|\sin\pi(a+\theta_j)|+\sum_{ j\in I'_1\cup  I'_2  , j\neq i }\ln |\sin\pi(a-\theta_j)|
+(2q_n-1)\ln2  $$
\begin{equation}\label{G36}
    =\Sigma_{+}+\Sigma_-+(2q_n-1)\ln2,  \;\;\;\;\;\;\;\;\;\;\;\;\;\;\;\;\;\;\;\;\; \;\;\;\;\;\;\;\;\;\;\;\;\;\;\;\;\;\;\;\;
\end{equation}
 where
 \begin{equation}\label{G37}
   \Sigma_{+}=\sum_{j\in I'_1\cup  I'_2  , j\neq i }\ln |\sin\pi(a+  \theta+j  \alpha)|,
 \end{equation}
 and
 \begin{equation}\label{G38}
     \Sigma_-=\sum_{j\in I'_1\cup  I'_2  , j\neq i }\ln |\sin\pi ( a-\theta   -j\alpha)|.
 \end{equation}
 Write $  \Sigma_{+}$ as the following form:
 \begin{equation}\label{G39}
   \Sigma_{+}=\sum_{j\in I'_1, j\neq i }\ln |\sin\pi(a+  \theta+j  \alpha)|+\sum_{j\in I'_2  }\ln |\sin\pi(a+  \theta+j  \alpha)|.
 \end{equation}
  We will estimate $\sum_{j\in I'_1, j\neq i }\ln |\sin\pi(a+  \theta+j  \alpha)|$ and $\sum_{j\in I'_2  }\ln |\sin\pi(a+  \theta+j  \alpha)|$
 respectively.
 \par
 On the one hand,
  $$\sum_{j\in I'_1, j\neq i }\ln |\sin\pi(a+  \theta+j  \alpha)| \;\;\;\;\; \;\;\;\;\; \;\;\;\;\; \;\;\;\;\; \;\;\;\;\;
   \;\;\;\;\; \;\;\;\;\; \;\;\;\;\; \;\;\;\;\; \;\;\;\;\; \;\;\;\;\; \;\;\;\;\;$$
   $$=  \sum_{j\in I'_1  }\ln |\sin\pi(a+  \theta+j  \alpha)|-\ln |\sin\pi(a+  \theta+i  \alpha)| \;\;\;\;\; \;\;\;\;\; \;\;\;\;\; \;\;\;\;
   $$
   $$  =  \sum_{ k=1  }^{q_n}\ln |\sin\pi(x +  k \alpha)|-\ln |\sin\pi(a+  \theta+i  \alpha)|  \;\;\;\;\; \;\;\;\;\; \;\;\;
    \;\;\;\;\; \;\;\;\;\;\;\; $$
    $$  \;\;\;\;\; \;\;\;\; =  \sum_{ k=1,k\neq k_0  }^{q_n}\ln |\sin\pi(x +  k \alpha)|+ \ln |\sin\pi(x +  k_0 \alpha)|-\ln |\sin\pi(a+  \theta+i  \alpha)|, $$
  where $x =a+\theta-([\frac{2}{3}q_n]+1)\alpha $ and $k_0$ satisfies
  $|\sin  \pi(x +  k_0   \alpha   )|=
\min_{1\le k\le q_n}|\sin  \pi(x +  k \alpha )|$.
In Lemma \ref{Le32}, let $m_k=0$, $k=1,2, \cdots q_n$, by the second equality of    (\ref{G34}),
 one has
 \begin{equation*}
 \sum_{ k=1,k\neq k_0  }^{q_n}\ln |\sin\pi(x +  k \alpha)|\leq  -(q_n-1)\ln 2+
C\ln q_n.
 \end{equation*}
Since   $\ln |\sin\pi(x +  k_0 \alpha)|\leq \ln |\sin\pi(a+  \theta+i  \alpha)| $ (by the minimality of $k_0$), we have
   \begin{equation}\label{G310}
   \sum_{j\in I'_1, j\neq i }\ln |\sin\pi(a+  \theta+j  \alpha)| \leq  -(q_n-1)\ln 2+
C\ln q_n.
   \end{equation}
   On the other hand,
   $$\sum_{j\in I'_2 }\ln |\sin\pi(a+  \theta+j  \alpha)|\;\;\;\;\;\;\;\;\;\;\;\;\;\;\;\;\;\;\;\;\;\;\;\;\;\;\;\;\;\;\;\;\;\;\;\;\;\;\;\;\;\;\;\;\;\;\;\;\;\;\;\;\;\;\;\;\;\;\;\;\;\;\;\;\;\;\;\;\;\; $$
   $$  =  \sum_{ k=1  }^{q_n}\ln |\sin\pi(x + ( k +m_k)\alpha)| \;\;\;\;\;\;\;\;\;\;\;\;\;\;\;\;\;\;\;\;\;\;\;\;\;\;\;\;\;\;\;\;\;\;\;\;\;\;\;\;\;\;\;\;\;\;  $$
    $$ =  \sum_{ k=1,k\neq k_0  }^{q_n}\ln |\sin\pi(x + ( k +m_k) \alpha)|+ \ln |\sin\pi(x +  ( k_0 +m_{k_0}) \alpha)|, $$
  where $x =a+\theta+(-[\frac{2}{3}q_n]+q_n-1)\alpha $, $m_k=0$ for $1\leq k\leq 2[\frac{2}{3}q_n]-q_n-1$ and
   $m_k=\ell -1$ for $ 2[\frac{2}{3}q_n]-q_n \leq k\leq  q_n $,  and $k_0$ satisfies
  $|\sin  \pi(x +  (k_0  +m_{k_0} \alpha   )|=
\min_{1\le k\le q_n}|\sin  \pi(x + ( k+m_k) \alpha )|$.
By the second equality of    (\ref{G34}) again,
 one has
 \begin{equation*}
 \sum_{ k=1,k\neq k_0  }^{q_n}\ln |\sin\pi(x +  (k+m_k) \alpha)|\leq  -(q_n-1)\ln 2+
C\ln q_n.
 \end{equation*}
In addition  $\ln |\sin\pi(x + ( k_0 +m_{k_0}\alpha)|\leq  0 $,  one has
  \begin{equation}\label{G311}
   \sum_{j\in I'_2 }\ln |\sin\pi(a+  \theta+j  \alpha)| \leq  -(q_n-1)\ln 2+
C\ln q_n.
  \end{equation}
  Putting  (\ref{G39}), $ (\ref{G310})$  and $ (\ref{G311})$ together, we have
  \begin{equation}\label{G312}
  \Sigma_{+} \leq  -2q_n\ln 2+
C\ln q_n.
  \end{equation}
   We are now  in the position to  estimate  $\Sigma_-$. In order to avoid repetition,  we omit some details. Similarly,
  $\Sigma_-$ consists of $2$ terms
of  the form  as $(\ref{G34})$,  plus two terms of the form $\min_{k=1,\ldots,q_{n}}\ln
|\sin  \pi(x+ (k+m_kq_n)\alpha   )|,$
where $m_k\in\{0,  (\ell-1)\},\;
k=1,\cdots,q_{n},$ minus $\ln
|\sin \pi (a -\theta_i )| $.  Following the estimate of $\Sigma_+ $,
 \begin{equation}\label{G313}
  \Sigma_{-} \leq  -2q_n\ln 2+
C\ln q_n.
  \end{equation}
Putting  $(\ref{G312})$  and  $(\ref{G313})$ into $(\ref{G36})$, we obtain
\begin{equation}\label{G314}
\sum_{{j\in I_1\cup I_2} {j\not=i}}\ln |\cos 2\pi a-\cos
2\pi\theta_j|\le  -2q_{n} \ln 2 +C\ln q_n.
\end{equation}
The estimate of $\sum _{ j\in I'_1\cup  I'_2  , j\neq i }\ln|\cos2\pi\theta_i  -\cos2\pi \theta_j|$ require a bit more work.
\par
It is easy to see that
 $$ \sum _{j \in I'_1  \cup I'_2,j\neq i}\ln|\cos2\pi \theta_i-\cos2\pi \theta_j| \;\;\;\;\;\;\;\;\;\;\;\;\;\;\;\;\;\;\;\;\;\;\;\;\;\;\;\;\;\;\;\;\;\;\;\;\;\;\;\;\;\;\;\;\;\;\;\;\;\;\;\;\;\;\;\;$$
\begin{equation}\label{G315}
    =\Sigma'_{+}+\Sigma'_-+(2q_n-1)\ln2,  \;\;\;\;\;\;\;\;\;\;\;\;\;\;\;\;\;\;\;\;\; \;\;\;\;\;\;\;\;\;\;\;\;\;\;\;\;\;\;
\end{equation}
where
\begin{equation}\label{G316}
   \Sigma'_{+}=\sum_{j \in I_1  \cup I_2,j\neq i}\ln |\sin\pi(2\theta+ (i+j) \alpha)|,
 \end{equation}
 and
 \begin{equation}\label{G317}
     \Sigma'_-=\sum_{j \in I_1  \cup I_2,j\neq i}\ln |\sin\pi( i-j)\alpha|.
 \end{equation}

 \par
Firstly,  we estimate $\Sigma'_+$. Similarly,   $\Sigma'_+$ consists of $2$ terms
of  the form  as $(\ref{G34})$,  plus two terms of the form $\min_{k=1,\ldots,q_{n}}\ln
|\sin  \pi(x+ (k+m_kq_n)\alpha   )|,$
where $m_k\in\{0,  (\ell-1)\},\;
k=1,\cdots,q_{n},$ minus $\ln
|\sin 2\pi (\theta+i\alpha) | $.
\par
   Following the  above arguments  and using  the first inequality  of  (\ref{G34}), we obtain
\begin{equation}\label{G318}
  \Sigma'_{+} >  -2q_n\ln 2-C\ln q_n+ 2\min_{j,i\in I_1\cup I_2} \ln |\sin
 \pi(2\theta+ (j+i )\alpha) |.
  \end{equation}
Thus  it is enough to    estimate the  last term in (\ref{G318}).  By  the hypothesis  $\theta\notin  \mathscr{R}$,
 one has
 \begin{equation}\label{G319}
  \min_{j,i\in [-2q_n,2q_n-1]} |\sin
 \pi(2\theta+ (j+i )\alpha) |>  \frac{1}{16 q_n^2}  \text{ for large n}.
 \end{equation}
 If $k\in I_2$, let $\ell_k=\ell-1$ and $ k'=k- \ell_kq_n$; if $k\in I_1$, let $\ell_k=0$ and $ k'=k$, then $i',j'\in  [-2q_n,2q_n-1]$.
 Recall that  $\Delta_n\leq\frac{1}{q_{n+1}}$.  It is easy to verify   $|\ell_k \Delta_n|<\frac{1}{ q_n^5}$ for $n$ large enough
 since $\beta(\alpha)>0$. Combining with  $(\ref{G319})$, we have for any $i,j \in I_1 \cup I_2$,
 \begin{equation*}
    |\sin  \pi(2\theta+ (j+i )\alpha) |\;\;\;\;\;\;\;\;\;\;\;\;\;\;\;\;\;\;\;\;\;\;\;\;\;\;\;\;\;\;\;\;\;\;\;\;\;\;\;\;\;\;\;\;\;\;\;\;\;
    \;\;\;\;\;\;\;\;\;\;\;\;\;\;\;\;\;\;\;\;\;\;\;\;\;\;\;\;\;\;\;\;\;\;
 \end{equation*}
 \begin{equation*}
    =|\sin\pi(2\theta+ (j'+i')\alpha )\cos \pi (\ell_i+\ell_j)\Delta_n
\pm \cos  \pi(2\theta+ (j'+i')\alpha ) \sin\pi (\ell_i+\ell_j)\Delta_n|
 \end{equation*}
\begin{equation}\label{G320}
 >\frac{1}{100q^2_n} \;\;\;\;\;\;\;\;\;\;\;\;\;\;\;\;\;\;\;\;\;\;\;\;\;\;\;\;\;\;\;\;\;\;\;\;\;\;\;\;\;\;\;\;\;\;\;\;\;
    \;\;\;\;\;\;\;\;\;\;\;\;\;\;\;\;\;\;\;\;\;\;\;\;\;\;\;\;\;\;\;\;\;\;\;\;\;\;\;\;\;\;\;\;
\end{equation}
(the $\pm$ depending on   the sign of $q_n \alpha-p_n$).
\par
Thus, by (\ref{G318}) and (\ref{G320}),
\begin{equation}\label{G321}
\Sigma'_+ >-2q_n\ln2-
C\ln q_n.
\end{equation}

Similarly, $\Sigma'_-$ consists of $2$ terms of the form   as  $ (\ref{G34}) $ plus  the  minimum term
(   because $\min_{j \in I'_1 }|\sin\pi( i-j)\alpha|=0$, then  $\sum  _{j \in I'_1,   j\neq i}\ln |\sin\pi( i-j)\alpha|$ is exactly of  the form   $ (\ref{G34}) $ ).
It follows that
\begin{equation}\label{G322}
  \Sigma'_{-} >  -2q_n\ln 2-C\ln q_n+  \min_{ j\in I_1\cup I_2, j\neq i} \ln |\sin
 \pi( (j-i )\alpha) |.
  \end{equation}
 We   are now in the position to estimate the last  term in (\ref{G322}).
Notice that for any $ i\in I_1\cup I_2$,  there is only one  $  \tilde{i}\in I_1\cup I_2$
such that $|i-\tilde{i}|=q_n$ or $\ell q_n$.
It is easy to check
\begin{equation}\label{G323}
   \ln|\sin\pi(i-\tilde{i})\alpha| \geq\min\{\ln|\sin\pi q_n\alpha|, \ln|\sin \pi \ell q_n\alpha|\}>-\ln q_{n+1}-C,
\end{equation}
since    $\Delta_n\geq \frac{1}{2q_{n+1}} $.
 If $j\neq i,\tilde{i}$ and $j\in I_1\cup I_2$, then $j-i=r+m'_jq_n$ with $1\leq |r|<q_n$ and $|m'_j|\leq \ell +2$. Thus
 by $(\ref{G31})$ and $(\ref{G32})$,
 $$||r\alpha||_{\mathbb{R}/ \mathbb{Z}}\geq \Delta_{n-1} \geq \frac{1}{2q_n}$$
 and
 \begin{eqnarray}
 \nonumber
    \min_{ {j\in I_1\cup I_2} {j\not=i,\tilde{i}}}\ln|\sin\pi(j-i)\alpha| & > & \ln (||r\alpha||_{\mathbb{R}/ \mathbb{Z}}- (\ell +2) \Delta_n )-C \\
   & >&  -\ln q_n -C, \label{G324}
 \end{eqnarray}
since $ (\ell +2) \Delta_n  <\frac{1}{10q_n}$ for $n$ large enough.

 By   $( \ref{G323})$ and  $ (\ref{G324})$, one has
 \begin{equation}\label{G325}
    \min_{ {j\in I_1\cup I_2} {j\not=i}}\ln|\sin\pi(j-i)\alpha| >  -\ln q_{n+1}-C\ln q_n.
 \end{equation}
By  the definition   $ \beta =\limsup_{n\rightarrow\infty}\frac{\ln q_{n+1}}{q_n}$,  (\ref{G322}) becomes
\begin{eqnarray}\label{G326}
\nonumber
    \Sigma'_-&>&-2q_n\ln2-\ln q_{n+1}-
 C\ln q_n \\
&>&-2q_n\ln2-(\beta + \varepsilon  )q_n-C\ln q_n,
\end{eqnarray}
for large $n$.
\par

By (\ref{G315}), (\ref{G321}) and (\ref{G326}),
 \begin{equation}\label{G327}
 \sum _{j \in I'_1  \cup I'_2,j\neq i}\ln|\cos2\pi \theta_i-\cos2\pi \theta_j|   > -2q_n\ln2-(\beta + \varepsilon  )q_n-C\ln q_n.
 \end{equation}
Together  with  (\ref {G314}),
we obtain
 $$\sum _{ j\in I'_1\cup  I'_2  , j\neq i }\ln|\cos2\pi a-\cos2\pi \theta_j|-  \ln|\cos2\pi\theta_i  -\cos2\pi \theta_j|< (\beta + \varepsilon  )q_n +C\ln q_n. $$
 This implies
$$\max_{ x\in[-1,1]}\max_{i=1,\cdots,k+1}\prod_{ j=1 , j\neq i }^{k+1}\frac{|x-\cos2\pi\theta_j|}
        {|\cos2\pi\theta_i-\cos2\pi\theta_j|}<e^{(2q_n-1)(\frac{\beta}{2}+\varepsilon ) }$$
         for large enough $ n$.$\qed$
         \par
        In Lemma $\ref{Le23}$, let $k=2q_n-1$,  $\epsilon_1=\frac{\beta}{2}+\varepsilon $ and $ \epsilon= \frac{\beta}{2}+2\varepsilon $.
        Clearly, $\epsilon_1< \epsilon$. Thus for any $\varepsilon >0$,  there exists some $j_0\in I_1\cup I_2$ such that $ \theta_{j_0}\notin A_{2q_n-1, \ln\lambda- \frac{\beta}{2}-2\varepsilon }$  for  $n$   large enough.

         \par\textbf{Step 2:} We will  show that $ \theta_j\in A_{2q_n-1, \ln\lambda- \frac{\beta}{2}-2\varepsilon } $ for all $j\in I_1$.
\begin{lemma}\label{Le34}
  $ \forall \varepsilon>0$, suppose  $ k\in [-2q_n,2q_n ]$ and  $d= dist(k,\{ m q_n\}_{m\geq0})\geq  \frac{q_n}{4}$,
then for sufficiently large
$n $
\begin{equation}\label{G328}
    |\phi(k)|<\exp(-(L-\varepsilon)d).
\end{equation}
\end{lemma}
\textbf{Proof: } We will use  block resolvent expansion   to prove this lemma. For any $ \varepsilon_0>0$, by hypothesis
$ k\in [-2q_n,2q_n ]$,  there exists   some $m\in \{-2,-1,0,1\}$ such that  $m q_n\leq k\leq (m+1)q_n$. $\forall y \in [m q_n+\varepsilon_0q_{n }+1,(m+1) q_n-
\varepsilon_0 q_{n }-1]$,  apply   Theorem $\ref{Th24}$ with  $\varepsilon=\varepsilon_0$, then    $sq_{n-1}\geq \frac{1}{2} dist(y,\{ m q_n\}_{m\geq0}) \geq \frac{\varepsilon_0q_n}{2}$
and
 \begin{equation*}
  \ln\lambda+9\ln (s q_{n-1}/q_n)/q_{n-1}-\varepsilon_0\geq \ln\lambda +9\frac{\ln(\varepsilon_0/2)}{q_{n-1}}-\varepsilon_0\geq  \ln\lambda-2\varepsilon_0,
 \end{equation*}
 for large $n$. Moreover,
 there exists an interval $ I(y)=[x_1,x_2]\subset
[ (m-1 ) q_n,(m+2)q_n]$
such that $y\in I(y)$ and
\begin{equation}\label{G329}
    \text{dist}(y,\partial I(y))\geq  \frac{1}{5} |I(y)|=\frac{2sq_{n-1}-1}{5} >\frac{q_{n-1}}{3}
\end{equation}
and
\begin{equation}\label{G330}
  |G_{I(y)}(y,x_i)| < e^{-(L- 2\varepsilon_0)|y-x_i|},\;i=1,2,
\end{equation}
where
 $ \partial I(y)$ is the boundary of the interval $I(y)$, i.e.,$\{x_1,x_2\}$, and recall that $ |I(y)|$ is the  number of $I(y) $, i.e., $ |I(y)|=x_2-x_1+1$.
   For $z  \in  \partial I(y)$,  let
  $z' $ be the neighbor of $z$, (i.e., $|z-z'|=1$) not belonging to $I(y)$.
\par
If $x_2+1<(m+1)q_n-\varepsilon_0q_n$ or  $x_1-1> m q_n+\varepsilon_0q_n$,
we can expand $\phi(x_2+1)$ or $\phi(x_1-1)$ as ($ \ref{G210}$). We can continue this process until we arrive to $z$
such that $z+1\geq(m+1)q_n-\varepsilon_0q_n$ or  $z-1\leq m  q_n+\varepsilon_0q_n$, or the iterating number reaches
$[\frac{3 d}{q_{n-1}}]$. Thus, by (\ref{G210})
\begin{equation}\label{G331}
   \phi(k)=\displaystyle\sum_{s ; z_{i+1}\in\partial I(z_i^\prime)}
G_{I(k)}(k,z_1) G_{I(z_1^\prime)}
(z_1^\prime,z_2)\cdots G_{I(z_s^\prime)}
(z_s^\prime,z_{s+1})\phi(z_{s+1}^\prime),
\end{equation}
where in each term of the summation one has
$m q_n+\varepsilon_0 q_{n }+1<z_i<(m+1) q_n-\varepsilon_0
q_{n }-1$, $i=1,\cdots,s,$ and
  either $z_{s+1} \notin [m q_n+\varepsilon_0 q_{n }+2,(m+1) q_n-
\varepsilon_0q_{n }-2]$, $s+1 < [\frac{3 d}{q_{n-1 }}]$; or
$s+1= [\frac{3 d}{q_{n-1}}]$.
\par
 If $z_{s+1} \notin [m q_n+\varepsilon_0q_{n }+2,(m+1) q_n-
\varepsilon_0 q_{n }-2]$, $s+1 < [\frac{3d}{q_{n-1 }}]$,
 by  ($\ref{G330}$),
\begin{eqnarray}
\nonumber
 & | G_{I(k)}(k,z_1) G_{I(z_1^\prime)}
(z_1^\prime,z_2)\cdots G_{I(z_s^\prime)}
(z_s^\prime,z_{s+1})\phi(z_{s+1}^\prime)|\;\;\;\;\;\;\;\;\;\;\;\;\;\;\;\;\;\;\;\;\;\;\;\;\;\;\;\;\;\;\;\;\;\;\;\;\;\\
\nonumber
&< e^{-(\ln\lambda-2\varepsilon_0)(|k-z_1|+\sum_{i=1}^{s}|z_i^\prime-z_{i+1}|)}
q_n ^C\;\;\;\;\;\;\;\;\;\;\;\;\;\;\;\;\;\;\;\;\;\;\;\;\;\;\;\;\;\;\;\;\; \;\;\;\;\\
&< e^{-(\ln\lambda-2\varepsilon_0)(|k-z_{s+1}|-(s+1))}q_n ^C <  e^{-(\ln\lambda-2\varepsilon_0)(d- \varepsilon_0q_{n } -4-\frac{ 3d}{q_{n-1}})}q_n^C,\;\;\; \label{G332}
\end{eqnarray}
since  $|\phi(z_{s+1}^\prime)|\leq(1+|z_{s+1}^\prime|)^C\leq q_n^C $.
If $s+1= [\frac{3 d}{q_{n-1}}] ,$
using   ($\ref{G329}$) and ($\ref{G330}$), we obtain
\begin{equation}\label{G333}
     | G_{I(k)}(k,z_1) G_{I(z_1^\prime)}
(z_1^\prime,z_2)\cdots G_{I(z_s^\prime)}
(z_s^\prime,z_{s+1})\phi(z_{s+1}^\prime)|<  e^{-(\ln\lambda-2\varepsilon_0) {\frac{q_{n-1}}{3}} [\frac{3 d}{q_{n-1}}]} q_n^C.
\end{equation}

Finally,   notice that the total number of terms in ($ \ref{G331}$)
is  at most  $2^{[\frac{3 d}{q_{n-1}}]}$ and $d\geq\frac{q_n}{4}$. Combining with
($\ref{G332}$) and
($\ref{G333}$),  we obtain
$$
|\phi(k)|<   e^{-(\ln\lambda-3\varepsilon_0-8\varepsilon_0\ln\lambda) d }
$$
for large $n$. By the arbitrariness of $\varepsilon_0 $,   we complete the proof of the lemma.
\begin{remark}\label{Remark410}
Under the hypothesis  of Lemma $\ref{Le34}$, Avila and Jitomirskaya  only   prove  that  $|\phi(k)|<\exp(-(\ln \lambda-\varepsilon)\frac{d}{2})$.
We give the refined version.
\end{remark}

\begin{theorem}\label{Th36}
 $ \forall \varepsilon>0$ and  for any $b\in[-\frac{5}{3} q_n,-\frac{1}{3}q_n]\cap \mathbb{Z}$,
we have $\theta+(b+q_n-1)\alpha\in A_{2q_n-1,2\ln\lambda/3+\varepsilon} $  if    $n$ is  large enough,  i.e., for all $j\in I_1$,
$\theta_j\in  A_{2q_n-1,2\ln\lambda/3+\varepsilon} $.
\end{theorem}
\textbf{Proof:} Let $b_1=b-1$ and $b_2=b+2q_n-1$.
 For any $\varepsilon_0>0 $, applying Lemma $\ref{Le34}$ (let $\varepsilon=\varepsilon_0$), one obtains  that   for $i=1,2$,
$$|\phi(b_i)|\leq \left\{
                \begin{array}{ll}
                  e^{-(\ln\lambda-\varepsilon_0)(2q_n+ b)}, &  {   -\frac{5q_n}{3}\leq b\leq -\frac{3q_n}{2};} \\
                  e^{-(\ln\lambda-\varepsilon_0)| q_n+b |}, &   {-\frac{3q_n}{2}< b< -\frac{ q_n}{2}\;and \;|b+q_n|>\frac{1}{4}q_n;} \\
                  e^{ (\ln\lambda-\varepsilon_0)  b }, &  {   -\frac{ q_n}{2}\leq b \leq -\frac{ q_n}{3}.}
                \end{array}
              \right.
$$
In  ($  \ref{G210}$), let  $I=[b,b+2q_n-2]$ and $x=0$, we get for $n$ large enough,
$$\max(|G_I(0,b)|,|G_I(0,b+2q_n-2)|)\geq
\left\{
  \begin{array}{ll}
     e^{ (\ln\lambda-2\varepsilon_0)(2q_n+ b)}, &  {   -\frac{5q_n}{3}\leq b\leq -\frac{3q_n}{2};} \\
    e^{ (\ln\lambda-2\varepsilon_0)| q_n+b |}, &  {  -\frac{3q_n}{2}< b< -\frac{ q_n}{2}\;and \;|b+q_n|>\frac{1}{4}q_n;} \\
    e^{ -(\ln\lambda-2\varepsilon_0)  b }, &  {   -\frac{ q_n}{2}\leq b\leq -\frac{ q_n}{3};} \\
     e^{-\varepsilon_0 q_n}, &  {   |b+q_n|\leq\frac{1}{4}q_n,}
  \end{array}
\right.
$$
since  $\phi(0)=1$ and $|\phi(k)|\leq (1+|k|)^C $.
\par
  Let $\varepsilon=\varepsilon_0$  in ($\ref{G25}$), and   let $ I=[b,b+2q_n-2]$, $y=0$, $k=2q_n-1$  in ($\ref{G26}$) and ($\ref{G27}$).  After
  careful computation, we obtain

   $|  Q_{2q_n-1}(\cos2\pi (\theta+(b+q_n-1)\alpha)|$

\par
$\;\;\;\;\;\;\;  =|P_{2q_n-1}(\theta+b\alpha)|$
\par
$ \;\;\;\;\;\;\; \leq\min\{|G_I(0,b)|^{-1}e^{(\ln\lambda+\varepsilon_0 )(b+2q_n-2)}, |G_I(0,b+2q_n-2)|^{-1}e^{-(\ln\lambda+\varepsilon_0)b}\} $
\par
$\;\;\;\;\;\;\; \leq e^{(2q_n-1)(2\ln\lambda/3+8\varepsilon_0)}.$

 By the arbitrariness of $\varepsilon_0 $, we finish the proof.$\qed$
 \par
 Since  $\ln\lambda> {\frac{3\beta}{2}}$,   $ \frac{2\ln\lambda}{3}< \ln\lambda-\frac{\beta}{2} $.
In  Step 1   and Step 2 if   let $\varepsilon $   be  so small that   $ \frac{2\ln\lambda}{3}+\varepsilon<\ln\lambda-\frac{\beta}{2}-2\varepsilon$,
i.e., $\varepsilon<\frac{1}{9} (\ln\lambda-\frac{3}{2}\beta)$,
  we have  $ \theta_j\in A_{2q_n-1, \ln\lambda- \frac{\beta}{2}-2\varepsilon } $ for all $j\in I_1$. This implies
there exists some  $j_0\in I_2$ such that $\theta_{j_0}\notin A_{2q_n-1,\ln\lambda-\frac{\beta}{2}-2\varepsilon }$
  if $\varepsilon<\frac{1}{9} (\ln\lambda-\frac{3}{2}\beta)$.
\par
\textbf{Step 3:} Establish the regularity for $y$.
\begin{theorem}\label{Th37}
 For any
$ \varepsilon>0$  such that $t= (\ln\lambda-\frac{3\beta}{2}-\varepsilon)>0$,
  $y$ is $(t,2q_n-1 )$-regular  for large enough $n$.
\end{theorem}
\textbf{Proof:}  According to the previous two steps,  there exists some $\theta_{j_0}\notin A_{2q_n-1,\ln\lambda-\frac{\beta}{2}-2\varepsilon_0 }$
for $j_0\in I_2$ if  $\varepsilon_0<\frac{1}{12} (\ln\lambda-\frac{3}{2}\beta)$.
 Set $I=[j_0-q_n+1,j_0+q_n-1]=[x_1,x_2]$.   In ($\ref{G25}$), let $\varepsilon=\varepsilon_0$,   combining with  ($\ref{G26}$)  and  ($\ref{G27}$),
 it is easy to verify
$$|G_I(y,x_i)|<e^{(\ln\lambda+\varepsilon_0 )(2q_n-2-|y-x_i|)-2q_n(\ln\lambda -\frac{\beta}{2}-2\varepsilon_0 )}.$$
By a simple computation $|y-x_i|\geq (\frac{2}{3}-\frac{1}{q_n^{1/9}} )q_n$,  then
 $$|G_I(y,x_i)|< e^{- |y-x_i|( \ln\lambda-\frac{3\beta}{2}-12\varepsilon_0 )},$$
 for large enough $n$.
  This implies  $y$ is  $(\ln\lambda-\frac{3\beta}{2}-12\varepsilon_0 ,2q_n-1)$-regular if   $\varepsilon_0<\frac{1}{12} (\ln\lambda-\frac{3}{2}\beta)$.
  For any $\varepsilon>0$  such that $t= (\ln\lambda-\frac{3\beta}{2}-\varepsilon)>0 $,  select $\varepsilon_0$ small enough  so that $ \ln\lambda-\frac{3\beta}{2}-\varepsilon<\ln\lambda-\frac{3\beta}{2}-12\varepsilon_0$.
  Then $y$ is $(t,2q_n-1)$-regular  for $n$ large
enough.
\section{The proof of Theorem $\ref{Th21}$ }
Now that  the regularity for   $y$ is  established, we will  use  block resolvent expansion   again to prove Theorem $\ref{Th21}$.
\par
\textbf{Proof of     Theorem $\ref{Th21}$.}

  Give some $k$ with  $   k>q_n$ and $n$ large enough.
$\forall y \in [ q_n^{\frac{8}{9}}, 2k]$,   let $\varepsilon=\varepsilon_0$  in   Theorem $ \ref{Th24}$  and $ \ref{Th37}$,  then there exists an interval $ I(y)=[x_1,x_2]\subset
[ -4k,4k]$ with
 $y\in I(y)$ such that
\begin{eqnarray}\label{G41}
\nonumber
    dist(y,\partial I(y)) &>& \frac{ 1}{5} |I(y)| \geq\min{\{\frac{2sq_{n-1}-1}{5},\frac{2q_n-1}{5} \}}  \\
    &\geq&  \frac{1}{3}q_{n-1}
\end{eqnarray}
and
\begin{equation}\label{G42}
  |G_{I(y)}(y,x_i)| < e^{-(\ln\lambda- \frac{3}{2}\beta-\varepsilon_0)|y-x_i|},\;i=1,2.
\end{equation}
As in  the proof of Lemma $\ref{Le34}$,
 denote by $ \partial I(y)$    the boundary of the interval $I(y)$. For $z  \in  \partial I(y)$,
let $z' $ be the neighbor of $z$, (i.e., $|z-z'|=1$) not belonging to $I(y)$.
\par
If $x_2+1< 2k$ or  $x_1-1>  b_n=q_n^{\frac{8}{9}}$,
we can expand $\phi(x_2+1)$ or $\phi(x_1-1)$ as ($ \ref{G210}$). We can continue this process until we arrive to $z$
such that $z+1\geq 2k$ or  $z-1\leq  b_n$, or the iterating number reaches
$[\frac{3 k}{q_{n-1}}]$.
\par  By (\ref{G210}),
\begin{equation}\label{G43}
   \phi(k)=\displaystyle\sum_{s ; z_{i+1}\in\partial I(z_i^\prime)}
G_{I(k)}(k,z_1) G_{I(z_1^\prime)}
(z_1^\prime,z_2)\cdots G_{I(z_s^\prime)}
(z_s^\prime,z_{s+1})\phi(z_{s+1}^\prime),
\end{equation}
where in each term of the summation we have
$b_n+1<z_i<2k-1$, $i=1,\cdots,s,$ and
  either $z_{s+1} \notin [b_n+2, 2k-2]$, $s+1 < [\frac{3 k}{q_{n-1 }}]$; or
$s+1= [\frac{3 k}{q_{n-1}}]$.
\par
 If $z_{s+1} \notin [b_n+2, 2k-2]$, $s+1 < [\frac{3k}{q_{n-1 }}]$, by  ($\ref{G42}$), one has
\begin{eqnarray}\label{G44}
\nonumber
 & | G_{I(k)}(k,z_1) G_{I(z_1^\prime)}
(z_1^\prime,z_2)\cdots G_{I(z_s^\prime)}
(z_s^\prime,z_{s+1})\phi(z_{s+1}^\prime)|\;\;\;\;\;\;\;\;\;\;\;\;\;\;\;\;\;\;\;\;\;\;\;\;\;\;\;\;\;\;\;\;\;\;\\
\nonumber
&\le e^{-(\ln\lambda- \frac{3}{2}\beta-\varepsilon_0)(|k-z_1|+\sum_{i=1}^{s}|z_i^\prime-z_{i+1}|)}
k ^C\;\;\;\;\;\;\;\;\;\;\;\;\;\;\;\;\;\;\;\;\;\;\;\;\;\;\;\;\;\;\;\;\;\;\; \;\;\;\; \; \\
\nonumber
&\le e^{-(\ln\lambda- \frac{3}{2}\beta-\varepsilon_0 )(|k-z_{s+1}|-(s+1))}k ^C\;\;\;\;\;\;\;\;\;\;\;\;\;\;\;\;\;\;\;\;\;\;\;\;\;\;\;\;\;\;\;\;\;\;\;\;\;\;\;\;\;\;\;\;\;\; \\
   &\le \max\{e^{-(\ln\lambda-\frac{3}{2}\beta-\varepsilon_0 )(k- b_n -4-\frac{ 3k}{q_{n-1}})}k^C ,
e^{-(\ln\lambda-   \frac{3}{2}\beta-\varepsilon_0)(2k-k -4-\frac{ 3k}{q_{n-1}})}k^C \}.
\end{eqnarray}
If $s+1= [\frac{3 k}{q_{n-1}}] ,$
using  ($\ref{G41}$) and ($\ref{G42}$),     we obtain
\begin{equation}\label{G45}
     | G_{I(k)}(k,z_1) G_{I(z_1^\prime)}
(z_1^\prime,z_2)\cdots G_{I(z_s^\prime)}
(z_s^\prime,z_{s+1})\phi(z_{s+1}^\prime)|\le  e^{-(\ln\lambda-\frac{3}{2}\beta-\varepsilon_0) {\frac{q_{n-1}}{3}} [\frac{3 k}{q_{n-1}}]} k^C.
\end{equation}

Finally,   notice that the total number of terms in ($ \ref{G43}$)
is  at most  $2^{[\frac{3 k}{q_{n-1}}]}$. Combining with
($\ref{G44}$) and
($\ref{G45}$),  we obtain
\begin{equation}\label{G46}
|\phi(k)|\le   e^{-(\ln\lambda-\frac{3}{2}\beta-2\varepsilon_0-\varepsilon_0\ln\lambda) k }
\end{equation}
for large enough $n$ (or equivalently  large enough $k$ ).  By the arbitrariness of $\varepsilon_0 $, we have for any $\varepsilon>0$,
\begin{equation}\label{G47}
|\phi(k)|\le   e^{-(\ln\lambda-\frac{3}{2}\beta-\varepsilon) k }  \text{  for } k \text{  large enough}.
\end{equation}
\par
For $k<0$, the proof is similar.
Thus  for any $\varepsilon>0$,
\begin{equation}\label{G48}
|\phi(k)|\le   e^{-(\ln\lambda-\frac{3}{2}\beta-\varepsilon) |k| }  \text{  if } |k| \text{  is large enough}.
\end{equation}
We finish the proof of    Theorem $\ref{Th21}$.

\begin{corollary}
Suppose  $\lambda>e^{\frac{3}{2}\beta}$ and $\theta \notin \mathscr{R} $. If   a solution  $\Psi_E(k)$
  satisfies   $H_{\lambda,\alpha,\theta}\Psi_E=E\Psi_E$ with $\Psi_E(k)\leq(1+|k|)^C$ and $E\in \Sigma_{\lambda,\alpha}$, then the  following holds:
\begin{equation}\label{G49}
    \limsup_{|k|\rightarrow\infty} \frac{\ln(\Psi_E^2(k)+\Psi_E^2(k+1))}{2|k|}\leq-(\ln \lambda- 3\beta /2  ).
\end{equation}
In particular, for $\beta(\alpha)=0$
\begin{equation}\label{G410}
    \lim _{|k|\rightarrow\infty} \frac{\ln(\Psi_E^2(k)+\Psi_E^2(k+1))}{2|k|}=- \ln \lambda  .
\end{equation}

\end{corollary}
\textbf{Proof:}     If  $\beta(\alpha)>0$, $ \forall \varepsilon>0$, by    $(\ref{G48})$,
\begin{equation*}
  |\Psi_E(k)|<e^{(\ln \lambda  - 3\beta /2- \varepsilon )|k|} \text{  for } |k| \text{  large enough}.
\end{equation*}
This implies
\begin{equation}\label{G411}
    \limsup_{|k|\rightarrow\infty} \frac{\ln(\Psi_E^2(k)+\Psi_E^2(k+1))}{2|k|}\leq- (\ln \lambda- 3\beta /2  )  \text{ if }\beta>0.
\end{equation}
\par
If $\beta(\alpha)=0$, following $\cite{AJ1}$ or $\cite{AJ2}$,  $k$ is $(t,  \ell(k))$-regular for large $|k|$, with $t=
\ln \lambda  - \varepsilon $. By the method  of  block resolvent expansion as above, we can obtain
\begin{equation*}
  |\Psi_E(k)|<e^{-(\ln \lambda - \varepsilon )|k|} \text{ if } k \text{ is large enough},
\end{equation*}
thus
\begin{equation}\label{G412}
    \limsup_{|k|\rightarrow\infty} \frac{\ln(\Psi_E^2(k)+\Psi_E^2(k+1))}{2|k|}\leq- \ln \lambda  .
\end{equation}
By (\ref{G411}) and (\ref{G412}), we obtain  (\ref{G49}).
\par
 By Furman's  uniquely ergodic Theorem   (Corollary 2   in   \cite{FUR} )
 \begin{equation}\label{G413}
    \liminf_{|k|\rightarrow\infty} \frac{\ln(\Psi_E^2(k)+\Psi_E^2(k+1))}{2|k|}\geq- \ln \lambda  .
\end{equation}
The last two inequalities   imply $( \ref{G410})$.
\begin{remark}
In $\cite{J}$, Jitomirskaya proves  $(\ref{G410})$  for  $ \alpha\in DC$, we extend   his result to all $\alpha$ with $\beta(\alpha)=0$.
 \end{remark}

           \begin{center}
           
             \end{center}
  \end{document}